\documentclass[12pt,graphicx]{amsart}

\usepackage{amscd,amssymb,amsthm,latexsym,amsfonts,graphicx,amsmath}
\newtheorem{thm}{Theorem}[section]
\newtheorem{lem}[thm]{Lemma}

\newtheorem{rem}{{\sc Remark}}[section]
\newcommand{\bas}{\begin{eqnarray*}}
\newcommand{\eas}{\end{eqnarray*}}
\newcommand{\ba}{\begin{eqnarray}}
\newcommand{\ea}{\end{eqnarray}}


\begin{document}
\title{On Constructing Special Lagrangian Submanifolds by Gluing}
%\footnote { Research supported in part by the Institute for Pure and Applied Mathematics under NSF grant DMS-9810282}

\author{Sema Salur}
\thanks{Research supported in part by the Institute for Pure and Applied Mathematics under NSF grant DMS-9810282}
%Preliminary Version
\vspace{.1in}
\address{Department of Mathematics, Cornell University, Ithaca, NY 14850}
\email{salur@math.cornell.edu}
\maketitle
%\footnotesize
\scriptsize

{\bf Abstract.}  The purpose of this paper is to give an application of the gluing theorem for special Lagrangian submanifolds of a Calabi-Yau 3-fold. In [2], a gluing theorem is proved to smooth a codimension-two singularity of a particular special Lagrangian submanifold. In this paper we will show that this theorem can be applied to more general cases where different special Lagrangians are intersecting and gives a way of constructing new special Lagrangian submanifolds. As an example we will show that a smooth special Lagrangian submanifold can be obtained from five copies of ${\bf RP^3}$ intersecting pairwise in a quintic.

\vspace{.1in}

\normalsize

\section{Introduction}

\vspace{.1in}

Due to SYZ conjecture [4], one important problem in Mirror Symmetry is to find a suitable compactification of the moduli space of special Lagrangian submanifolds. In particular, one should understand the singularities of this moduli space. For this purpose, in [2], we showed that a particular special Lagrangian submanifold with an irreducible singularity can be a limit point in this space by proving the following theorem:

\vspace{.1in}

{\thm: [2] Given a connected immersed special Lagrangian submanifold $L^3$ of a Calabi-Yau manifold $X^6$ with a particular irreducible, orthogonal self intersection $K$ of codimension-two (singularity of type $z_1.\overline {z_2}=0$) it can be approximated by a sequence of smooth special Lagrangian submanifolds and therefore $L$ is a limit point in the moduli space.}

\vspace{.1in}

In this paper we will show that Theorem 1.1 can be applied to more general cases where different special Lagrangians are intersecting and hence we can construct new special Lagrangian submanifolds nearby by gluing. In particular we will modify Theorem 1.1 as follows:
 
\vspace{.1in}

{\thm: Given a singular special Lagrangian submanifold $L$ which is invariant under a geometric $\bf {Z}_m$ action and consists of pairwise, orthogonal and cyclic intersections of special Lagrangian submanifolds $L_1,...,L_m$ of a three dimensional Calabi-Yau manifold $X$, it can be approximated by a sequence of smooth special Lagrangian submanifolds obtained by smoothing each codimension-two intersections $K_i$ (singularity of type $z_1.\overline {z_2}=0$) }

\vspace{.1in}

\section{The Eigenvalue Estimate}

In this section we will explain how Theorem 1.1 can be generalized to pairwise intersections of several special Lagrangian submanifolds.

\vspace{.1in}

In what follows, $X$ will denote a 3-dimensional Calabi-Yau manifold with $m$ different special Lagrangians $L_1,...,L_m$ intersecting pairwise (locally $z_i\cdot\overline z_{i+1}=0$ for $z_i\in L_i$, $1\leq i\leq m-1$) and perpendicularly with respect to the induced metric along curves $K_i$ for each $i$. Moreover we will assume that $L_m$ intersects $L_1$ to complete the cycle and $\bf{Z}_m$ is acting isometrically on $X$ and the singular special Lagrangian $L=\cup L_i$ is invariant under this action.

%also as in [2] each $L_i$ intersects with $L_{i+1}$ perpendicularly along a codimension two intersection $K_i$ for each $i$.

\vspace{.1in}

In [2], for a given singular special Lagrangian submanifold $L$ and a gluing parameter $\delta$, we first construct an approximate special Lagrangian submanifold $H_\delta$ in an open ball $V$ around the singular set $K$ and use the Implicit Function Theorem to prove that there exists a {\em true} special Lagrangian submanifold nearby. In order to prove the existence of a {\em true} special Lagrangian we need to get a uniform estimate for the right inverse for the linearized operator $D_{\delta}$. Here we need to do the same for each intersection.

{\rem:} Note that since we assume $K$ is irreducible, Theorem 1.1 is not strong enough to be applied directly to the case where we glue two arbitrary special Lagrangian submanifolds. Even though we need this irreducibility condition only in proving the eigenvalue estimate for the linearized operator and that other parts of the proof do not require this condition it is a crucial restriction on the gluing model. Without this assumption Theorem 1.1 is not true. However, one can automatically obtain irreducibility requirement if singular special Lagrangian is replaced by a group of intersecting totally symmetric special Lagrangian submanifolds which complete a cycle. Totally symmetric means that $\bf {Z}_m$ acting isometrically on $X$ and $L=\cup L_i\subset X$ is invariant under this action. Completing a cycle means that if $L$ consists of $m$ special Lagrangian submanifolds $L_1,...,L_{m}$ then $L_i$ intersects $L_{i+1}$ along curves $K_i$ for $1\leq i\leq m-1$ and $L_{m}$ intersects $L_1$ along $K_{m}$.

\vspace{.1in}

Given the gluing parameters $\delta_1,...,\delta_m$ we can smooth the singularities of the form $z_1\cdot \overline z_{i+1}=0$ for each intersection inside open neighbourhoods $V_i$ around $K_i$ and construct approximate special Lagrangians $H^i_{\delta_i}$ which agree with $L_i$ and $L_{i+1}$ outside a tubular neighbourhood of their intersection $K_i$. As before one can show that for each $H^i$ the linearized operator for the special Lagrangian equation is also $\Psi_i \cdot \Delta^i_{\delta_i}$ for each intersection $K_i$ where $\Psi_i$ is a small function for small values of $\delta_i$. Therefore it is sufficient to check the invertibility of $\Delta^i_{\delta_i}$. Also note that $\bf {Z_m}$ is acting on $L$ isometrically and therefore it is sufficient to check the invertibility for only one intersection. 

\vspace{.1in} 

Next, we will modify the eigenvalue estimates for each Laplacian operator $\Delta^i_{\delta_i}$ which are needed in the gluing theorem [2], and in section 3 we will apply this to an example.

{\lem : There are constants $C_i>0$ $(1\leq i\leq m)$independent of the gluing parameters $\delta_i$, such that for $\delta_i$ sufficiently small, the first (nonzero) eigenvalues $\lambda_1(\Delta^i_{\delta_i})$ of $\Delta^i_{\delta_i}$ are bounded below by $C_i$.}

\vspace{.1in}

{\bf Proof:} As in [2], we prove it by contradiction. Note that we have $m$ Laplacian operators ${\Delta}^1_{\delta_1},...,{\Delta}^m_{\delta_m}$ for $H^1,...,H^m$ and since $\bf{Z_m}$ is acting isometrically we can assume that they are all equivalent. Therefore the analysis reduces to the case where two special Lagrangian submanifolds are intersecting as before,[2]. 

%We will start with $H^1$ for $L_1\cap L_2$ and its corresponding operator $\Delta^1_{\delta_1}$ and the analysis for other intersections will follow because of invariance under $\bf{Z_m}$ action. 

\vspace{.1in}

Suppose that the lemma is not true for $\Delta^1_{\delta_1}$ in $L_1\cup L_2$. Then we may assume that the first eigenvalue $\lambda_1(\Delta^1_{\delta_1})$ converges to zero as $\delta_1$ tends to zero. Since we have the equivalence coming from the isometric action we can drop the index 1 in $\delta_1$. Let $\phi_{\delta}$ be the eigenfunction of $\lambda_1(\Delta^1_{\delta})$ satisfying

\vspace{.1in}

$\displaystyle \int_{\overline {H_{\delta_i}}} |\phi_\delta|^2 =1$ and $\displaystyle \int_{\overline {H_{\delta_i}}} \phi_\delta =0$ and $\Delta_{\delta}\phi_\delta=\lambda_{1,\delta}\phi_\delta$ .

\vspace{.1in}

\noindent where $\lambda_{1,\delta}$ determines the dependence of the first eigenvalue on the gluing parameter $\delta$ and ${\overline {H_{\delta_i}}}$ is the connected union of smoothed approximate special Lagrangians $H^i_{\delta_i}$. 

%Here for $H^1_{\delta_1}$ we will repeat the same limiting argument in [2] for $\delta_1$ but assume that when $\delta_1\rightarrow 0$ the other gluing parameters $\delta_2,...,\delta_m$ are fixed. 

%Note that every time we change $\delta$, we change the induced metric on the approximated special Lagrangian and since the Laplacian operator depends on the metric, the eigenvalues of $\Delta_{\delta}$ depend on $\delta$. For simplicity we will use $\lambda_{\delta}$ for $\lambda_{1,\delta}$.

\vspace{.1in}

For small compact sets away from singularity, the $L^2_\delta$ norm is uniformly equivalent to the usual $L^2$ norm. On these compact sets there exists a subsequence of $\phi_n$ that converges smoothly to a limit $\Delta \phi_0=0$. Following the same argument for the sequence of compact sets, and passing to a diagonal subsequence, we obtain a nonzero eigenfunction $\phi_0$ as the limit defined in the complement of the singularity satisfying 

\vspace{.1in}

$\displaystyle \int |\phi_0|^2 =1$ and $\displaystyle \int \phi_0 =0$

\vspace{.1in}

We now explain why $\phi_0$ cannot be zero.  If $\phi_0=0$ then for very small $\delta$, $\phi_{\delta}$ will be very small everywhere (almost zero) which contradicts the fact that 

\vspace{.1in}
$||\phi_{\delta}||_{L^2}\leq ||\phi_{\delta}||_{L^\infty}$ and our initial assumption $||\phi_{\delta}||_{L^2}=1$.

\vspace{.1in}

So we have a nonzero function $\phi_0$ in the limit and since $\lambda_{\delta}\rightarrow 0$ we get $\Delta_0\phi_0=0$. On a compact manifold the only harmonic functions are constant functions. Therefore $\phi_0$ should be some nonzero constant. On one component $\phi_{\delta}$ will converge to a constant and on the other component it will converge to another constant. Since $\bf{Z_m}$ is acting isometrically on $L$ these two constants should be same and since $\displaystyle \int \phi_0 =0$ this is only possible if $\phi_{\delta}$ converges to zero. This contradicts the fact that $\phi_0$ is nonzero.

\vspace{.1in}

One other possibility is the case when the eigenfunctions get trapped in the neck region and as the gluing parameter $\delta$ goes to 0 they converge to maps which are identically zero everywhere but blow up at one point. Here there is no need to study the concentration problem in the neck area because the analysis follows exactly the same way for each intersection as before [2].

%Since we assume that $\displaystyle \int \phi_0 =0$ and $\displaystyle \int |\phi_0|^2 =1$ this is only possible if $\phi_{\delta}$ converges to 1 in one component and -1 on the other. Without loss of generality we can assume that $\phi_{\delta}\rightarrow 1$ in $L_1$ and $\phi_{\delta}\rightarrow -1$ in $L_2$. 

%Since we assume that the singularity is irreducible these two constants should be same and since $\displaystyle \int \phi_0 =0$ this is only possible if $\phi_{\delta}$ converges to zero. This contradicts the fact that $\phi_0$ is nonzero. 

%\vspace{.1in}

%The arguments follow exactly the same for $L_2\cup L_3$ and $\Delta^2_{\delta_2}$ and similarly for the other intersections because $\bf{Z_m}$ is acting isometrically on $L$ and $L$ is invariant under this action. In each case we get that the eigenfunctions should converge to the same nonzero constant in each component which contradicts the fact that $\displaystyle \int \phi_0 =0$. 

%Hence the constantsThis gives us eigenfunctions converging to -1 on $L_2$ and 1 on $L_1$. When we do the same for the intersection $L_5\cap L_1$ we get $\phi_{\delta}$ converges to 1 on $L_5$. But this contradicts the fact that $\displaystyle \int \phi_0 =0$. Hence the constants should be all zero. 

\vspace{.1in}

Hence we can modify Theorem 1.1 as follows:
 
\vspace{.1in}

{\thm: Given a singular special Lagrangian submanifold $L$ which is invariant under a geometric $\bf {Z}_m$ action and consists of pairwise, orthogonal and cyclic intersections of special Lagrangian submanifolds $L_1,...,L_m$ of a three dimensional Calabi-Yau manifold $X$ (as in figure1), it can be approximated by a sequence of smooth special Lagrangian submanifolds obtained by smoothing each codimension-two intersections $K_i$ (singularity of type $z_1.\overline {z_2}=0$) }

\begin{figure}[h]

%\hspace{2 in}\includegraphics{cokgen2}

\includegraphics{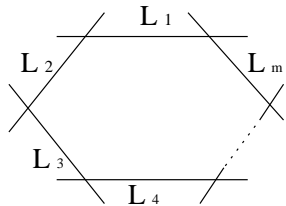}

\caption{\small $L=\cup L_i$}

%\label{sev-1}

\end{figure}

\vspace{.1in} 

{\rem :} Here we assumed that $L_1,...,L_m$ are intersecting orthogonally with respect to the induced metric but in the next example we have to verify that this is true. We will do this by averaging the metric with some finite group and making it invariant under this group action.

\section{The Example}

\vspace{.1in}

In this section we will apply Theorem 2.2 to five intersecting copies of ${\bf RP^3}$, [3], to obtain new special Lagrangian submanifolds in a quintic.

\vspace{.1in}

Let $X$ be a 3-dimensional Calabi-Yau manifold defined as a degree five hypersurface in ${\bf CP^4}$ given as follows:

\vspace{.1in} 

\hspace{.1in}  $X=\{ z_0^5+z_1^5+z_2^5+z_3^5+z_4^5=0\} \subset {\bf CP^4}$

\vspace{.1in} 

We will first write five different anti-holomorphic involutions $f_1$,..., $f_5$ on $X$. Then we will find the fixed point sets of these involutions and call them $F_1,...,F_5$. By a theorem of R. Bryant [1], $F_1$, ..., $F_5$ will be five different special Lagrangian submanifolds of the quintic $X$. Each of them can be visualized as the real part of $X$ and is diffeomorphic to ${\bf RP^3}$. In our example they also intersect pairwise as in figure 2. Moreover we will write a $\bf {Z}_5$ action on ${\bf CP^4}$ which acts isometrically on $X$. 

\begin{figure}[h]

\includegraphics{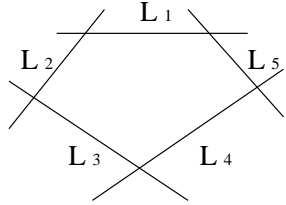}

\caption{\small Five different special Lagrangian submanifolds}

%\label{sev-1}

\end{figure}

\vspace{.1in} 

Next, we will justify this figure. Let $\xi=a+ib$,  $(a,b\in {\bf R})$ be the fifth root of unity. For $z_i=(x_i, y_i)$, let $L_1,L_2,L_3,L_4,L_5$ be defined as follows:

\vspace{.1in}

%$$ \{a\}$= fixed point set of$ \left\{ \begin{array}{ccc}

%z_0\rightarrow \overline z_0 \\

%z_1\rightarrow \overline z_1 \\

%z_2\rightarrow \overline z_2 \\

%z_3\rightarrow \overline z_3 \\

%z_4\rightarrow \overline z_4 \\

%\end{array} \right.  $$
\vspace{.1in} 

\noindent $$L_1= \left\{ \begin{array}{ccc}

$fixed point set of the involution$ \\

z_0\rightarrow \overline {z_0},$ $ z_1\rightarrow \xi\overline {z_1},$ $  z_2\rightarrow \xi\overline {z_2},$ $  z_3\rightarrow \overline {z_3},$ $  z_4\rightarrow \overline {z_4} \\
\end{array} \right.  $$

\vspace{.1in} 

=$\{x_0,x_1,x_2,x_3,x_4$ $|$ $x_0^5 + (x_1 +i\frac{1-a}{b}x_1)^5+ (x_2 +i\frac{1-a}{b}x_2)^5+ x_3^5 +x_4^5=0\}$

%\noindent $=\{ (x_0,0), (0,y_1), (0,y_2), (x_3,0), (x_4,0)| x_0^5+x_3^5+x_4^5=0, y_1^5+y_2^5=0\}$ 
\vspace{.2in}

$$L_2= \left\{ \begin{array}{ccc}

$fixed point set of the involution$ \\

z_0\rightarrow \overline {z_0},$ $  z_1\rightarrow \overline {z_1},$ $ z_2\rightarrow \xi\overline {z_2},$ $  z_3\rightarrow \xi\overline {z_3},$ $  z_4\rightarrow \overline {z_4} \\

\end{array} \right.  $$

\vspace{.1in} 

=$\{x_0,x_1,x_2,x_3,x_4$ $|$  $x_0^5 + x_1^5 + (x_2 +i\frac{1-a}{b}x_2)^5+ (x_3 +i\frac{1-a}{b}x_3)^5+x_4^5=0\}$

%\noindent $=\{ (x_0,0), (x_1,0), (0,y_2), (0,y_3), (x_4,0)| x_0^5+x_1^5+x_4^5=0, y_2^5+y_3^5=0\}$

\vspace{.2in}

\noindent $$L_3= \left\{ \begin{array}{ccc}

$fixed point set of the involution$ \\

z_0\rightarrow \overline {z_0},$ $  z_1\rightarrow \overline {z_1},$  $  z_2\rightarrow \overline z_2,$ $  z_3\rightarrow \xi\overline {z_3},$ $  z_4\rightarrow \xi\overline {z_4} \\

\end{array} \right.  $$

\vspace{.1in}

%\noindent $=\{ (x_0,0), (x_1,0), (x_2,0), (0,y_3), (0,y_4)| x_0^5+x_1^5+x_2^5=0, y_3^5+y_4^5=0\}$ 

=$\{x_0,x_1,x_2,x_3,x_4$ $|$ $x_0^5 + x_1^5 + x_2^5 + (x_3+i\frac{1-a}{b}x_3)^5+ (x_4 +i\frac{1-a}{b}x_4)^5=0\}$

\vspace{.2in}

\noindent $$L_4= \left\{ \begin{array}{ccc}

$fixed point set of the involution$ \\

z_0\rightarrow \xi\overline {z_0},$ $ z_1\rightarrow \overline {z_1},$ $  z_2\rightarrow \overline z_2,$ $  z_3\rightarrow \overline {z_3},$ $  z_4\rightarrow \xi\overline {z_4}\\

\end{array} \right.  $$

\vspace{.1in} 

=$\{x_0,x_1,x_2,x_3,x_4$ $|$ $(x_0 +i\frac{1-a}{b}x_0)^5+x_1^5+ x_2^5+x_3^5+ (x_4 +i\frac{1-a}{b}x_4)^5=0\}$ 

%\noindent $=\{ (0,y_0), (x_1,0), (x_2,0), (x_3,0), (0,y_4)| x_1^5+x_2^5+x_3^5=0, y_0^5+y_4^5=0\}$ 

\vspace{.2in} 

\noindent $$ L_5= \left\{ \begin{array}{ccc}
$fixed point set of the involution$ \\

 z_0\rightarrow \xi\overline {z_0},$ $  z_1\rightarrow \xi\overline {z_1},$ $  z_2\rightarrow \overline {z_2},$ $ z_3\rightarrow \overline {z_3},$ $  z_4\rightarrow \overline {z_4} \\

\end{array} \right.  $$

\vspace{.1in} 

=$\{x_0,x_1,x_2,x_3,x_4$ $|$ $(x_0 +i\frac{1-a}{b}x_0)^5+ (x_1 +i\frac{1-a}{b}x_1)^5+x_2^5+x_3^5+x_4^5=0\}$ 

%\noindent $=\{ (0,y_0), (0,y_1), (x_2,0), (x_3,0), (x_4,0)| x_2^5+x_3^5+x_4^5=0, y_0^5+y_1^5=0\}$ 

\vspace{.2in}

\noindent The sets of intersection are as follows:

\vspace{.1in} 
Since 
\vspace{.1in} 

$$ L_1\cap  L_2= \left\{ \begin{array}{ccc}

z_0^5+z_1^5+z_2^5+z_3^5+z_4^5=0 \\

 z_0=\overline {z_0},$    $ z_1=\xi\overline {z_1},$   $ z_2=\xi\overline {z_2},$   $
 z_3=\overline {z_3},$   $ z_4=\overline {z_4} \\

 z_0=\overline {z_0},$   $z_1=\overline {z_1},$  $ z_2=\xi\overline {z_2},$ $ z_3=\xi\overline {z_3},$ $ z_4= \overline {z_4} \\

\end{array} \right.  $$

\vspace{.1in} 

\noindent this implies 

\vspace{.1in} 

$K_1$=  $L_1\cap  L_2= \{ z_1=0,$ $z_3=0,$ $z_0^5+z_2^5+z_4^5=0 \}\cong S^1$,

\vspace{.1in} 

\noindent and similarly we get

\vspace{.1in} 
$K_2$=  $L_2\cap L_3= \{ z_2=0,$ $z_4=0,$ $z_0^5+z_1^5+z_3^5=0 \}\cong S^1$,

\vspace{.1in} 
$K_3$=  $L_3\cap  L_4= \{ z_0=0,$ $z_3=0,$ $z_1^5+z_2^5+z_4^5=0 \}\cong S^1$,

\vspace{.1in} 
$K_4$=  $L_4\cap  L_5= \{ z_1=0,$ $z_4=0,$ $z_0^5+z_2^5+z_3^5=0 \}\cong S^1$,

\vspace{.1in} 
$K_5$=  $L_5\cap  L_1= \{ z_0=0,$ $z_2=0,$ $z_1^5+z_3^5+z_4^5=0 \}\cong S^1$,

\vspace{.1in} 

\noindent and for the other pairs we get

\vspace{.1in} 

$L_1\cap  L_3= \{ z_1=0,$ $z_2=0,$ $z_3=0$, $z_4=0$, $z_5=0 \}=\emptyset$,  

\vspace{.1in}

\noindent and similarly 

\vspace{.1in}

$L_1\cap  L_4= \emptyset$, $L_2\cap  L_4= \emptyset$, $L_2\cap  L_5= \emptyset$, and $L_3\cap L_5= \emptyset$.

\vspace{.1in} 
Next, we will describe the $\bf {Z_5}$$=\{g_0,g_1,g_2,g_3,g_4\}$ action which keeps $L=\cup L_i$ invariant.

\vspace{.1in} 

For all $i=0,1,...,4$,  $g_i$ induces a map $\tilde {g_i}: {\bf CP^4}\rightarrow {\bf CP^4}$ defined as:

\vspace{.1in} 
 $\tilde {g_0}=$id 
\vspace{.1in} 

$\tilde {g_1}:(z_0,z_1,z_2,z_3,z_4)\rightarrow(z_4,z_0,z_1,z_2,z_3)$,
\vspace{.1in} 

$\tilde {g_2}:(z_4,z_0,z_1,z_2,z_3)\rightarrow(z_3,z_4,z_0,z_1,z_2)$
\vspace{.1in} 

$\tilde {g_3}:(z_3,z_4,z_0,z_1,z_2)\rightarrow(z_2,z_3,z_4,z_0,z_1)$
\vspace{.1in} 

$\tilde {g_4}:(z_2,z_3,z_4,z_0,z_1)\rightarrow(z_1,z_2,z_3,z_4,z_0)$ are cyclic permutations.
 
\vspace{.1in} 
Since the involutions $f_i: {\bf CP^4}\rightarrow {\bf CP^4}$ satisfy $f_i\circ\tilde {g_i}= \tilde {g_i}\circ f_{i+1}$ for all $i$ and the quintic $X$ is invariant under the maps $\tilde {g_i}$, this implies that $\tilde {g_i}:L_i=F_i\cap X \rightarrow L_{i+1}=F_{i+1}\cap X $ where $F_i$ are the fixed point sets of the involutions $f_i$. One can also easily show that $\tilde {g_i}$ will take the intersections $K_i=L_i\cap L_{i+1}$ to $K_{i+1}=L_{i+1}\cap L_{i+2} $. These will imply that $L=\cup L_i$ is invariant under the ${\bf Z_5}$ action. 

\vspace{.1in}

As we mentioned in Remark 2.2 we need to verify that these special Lagrangian submanifolds intersect orthogonally with respect to the induced metric. Let $\tilde G$ be the group generated by the group of anti-holomorphic involutions and ${\bf Z_5}$. By construction it is a finite group and we can average the ambient metric so that it is invariant under the group action generated by $\tilde G$. This invariance and representation theory will then imply that $L_i\cap L_{i+1}$ orthogonally for each $i$.

\vspace{.1in}

%Also we can assume that the metric on the quintic $X$ is ${\bf Z_5}$ invariant and we can conclude that ${\bf Z_5}$ is acting isometrically on $L$. 

Then applying Theorem 2.2 we can smooth the singularities and obtain a smooth special Lagrangian submanifold which agrees with $L_1,...,L_5$ outside the balls $V_1,...,V_5$. 

{\rem :} By taking involutions appropriately in other Calabi Yau manifolds one can construct different special Lagrangian submanifolds using the same gluing process. These examples will be discussed somewhere else.

%\vspace{.1in} 

%{\rem :} We also hope to study the effect of $T$-duality transformation on this new special Lagrangian submanifold and possible applications in an upcoming paper.

\vspace{.2in}

\small

{\em Acknowledgements.} 

%This work was done when the author was visiting M.I.T. during the spring of 1999. Many thanks to the mathematics department at M.I.T. for their hospitality and support during the course of this work. 

This work was completed when the author was attending to the Conformal Field Theory program at the Institute for Pure and Applied Mathematics during Fall 2001. Many thanks to the IPAM and the mathematics department at UCLA for their support and hospitality. The author is also grateful to Paul Seidel for suggesting the example which motivated this paper at the $8^{th}$ Gokova Geometry and Topology Conference, May 2001.

\vspace{.4in} 
  

\begin{thebibliography}{[FP]}



%\bibitem[13]{bf} McLean, R.C. {\em Deformations of calibrated submanifolds}, Comm. Anal. Geom. {\bf 6} (1998), 705-747 

\bibitem[1]{bf} Bryant, R.L. {\em Minimal Lagrangian submanifolds of K{\"a}hler-Einstein
manifolds}, Differential geometry and differential equations
(Shanghai, 1985), 1--12, Lecture Notes in Math. {\bf 1255}, Springer,
Berlin-New York, 1987.  MR 87a:53082

\bibitem[2]{bf} Salur, S. {\em A Gluing Theorem for Special Lagrangian Submanifolds}, math.DG/0108182

\bibitem[3]{b} Seidel, P. private communication.


\bibitem[4]{bf} Strominger, A., Yau, S.T. and Zaslow, E., {\em Mirror Symmetry is T-Duality}, Nucl. Phys. {\bf B479} (1996), 243-259



\end{thebibliography}
\end{document}